\documentclass[12pt]{article}
\newcommand{\editor}{Renata Del-Vecchio}
\newcommand{\ArticleType}{Research Article}
\newcommand{\RecievedDate}{Mar 22, 2026}
\newcommand{\RevisedDate}{Aug 13, 2026}
\newcommand{\AcceptedDate}{Aug 14, 2026}
\newcommand{\PublishedDate}{Aug 15, 2026}
\newcommand{\JournalIndex}{Volume 5 (2026), Pages 81--93}
\newcommand{\LastName}{Prodinger}
\newcommand{\JournalIndexShort}{\LastName/ American Journal of Combinatorics 5 (2026) 81--93}
\usepackage{fancyhdr}
\usepackage{calrsfs}
\usepackage{mathrsfs}

\usepackage{amsmath,amsthm,amssymb,mathtools,url,graphicx,pdfpages,tikz,rotating}
\usepackage[bookmarks=true,pdfborder={0 0 0}]{hyperref}
\usepackage[affil-it]{authblk}
\usepackage[left=1in,right=1in,top=1.2in, bottom=1in]{geometry}

\theoremstyle{definition}

\numberwithin{equation}{section}

\begin{document}
\setcounter{page}{81}
\noindent {\color{teal}\bf\large American Journal of Combinatorics} \hfill \ArticleType\\
\JournalIndex

\title{Ordered trees with distinguished children}

\author{Helmut Prodinger}
\affil{\normalsize\rm (Communicated by \editor)
\vspace*{-24pt}}
\date{}

{\let\newpage\relax\maketitle}

\begin{abstract}
A new tree model is introduced based on ordered trees, by distinguishing exactly one child of each node that \emph{has} children.
The basic enumeration leads to a cubic equation of the generating function. The extraction of its coefficients can be done
using the Lagrange inversion formula. Various parameters that are commonly studied for ordered trees can also be addressed here, like
degree of the root, number of leaves, number of old leaves, height, height of leftmost leaf, and pathlength. We go through these instances and leave
further parameters to later research, by either the author or some readers. Dealing with cubic equations is essential.
Finally,  ordered trees are replaced by marked ordered trees; they are then combined with the concept of distinguished children. Only the basic 
enumeration is provided at that stage.
\end{abstract}

\renewcommand{\thefootnote}{\fnsymbol{footnote}} 
\footnotetext{\hspace*{-22pt} MSC2020: 05A15;
Received \RecievedDate; Revised \RevisedDate; Accepted \AcceptedDate; Published \PublishedDate \\
\copyright\, The author(s). Released under the CC BY 4.0 International License}
\renewcommand{\thefootnote}{\arabic{footnote}}

\section{Old and young leaves and the present paper}

Plane trees (ordered trees) have a root and a sequence of $k\ge0$ subtrees, which are (recursively) plane trees.
The enumeration is with the respect to the number of nodes or (equivalently) the number of edges and leads to
a prominent structure related to \emph{Catalan} numbers. Leaves are nodes that have no successors. By agreement,
the tree just consisting of one node is not considered as a leaf. The enumeration with respect to nodes and leaves leads
to \emph{Narayana} numbers. A further refinement is due to Chen, Deutsch, and Elizalde~\cite{CDE}: Leaves that are the
leftmost child of its parent is called \emph{old}, all the other leaves are called \emph{young}.
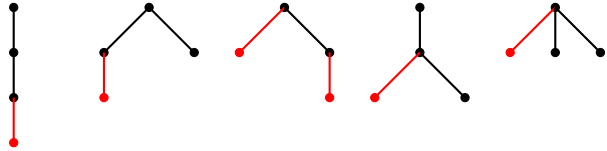
\begin{figure}[h]
\begin{center}
	\begin{tikzpicture}[scale=0.6]

		\draw[black,fill=black,xshift=3cm] (0,0) circle (.5ex);
		\draw[black,fill=black,xshift=3cm] (0,-1) circle (.5ex);
		\draw[black,fill=black,xshift=3cm] (0,-2) circle (.5ex);
		
		\draw [ thick, xshift=3cm] (0,0) -- (0,-1) ;
		\draw [thick,xshift=3cm] (0,-1) -- (0,-2) ;
		\draw [thick,red,xshift=3cm] (0,-2) -- (0,-3) ;
		\draw[red,fill=red,xshift=3cm] (0,-3) circle (.5ex);
		
		\draw[black,fill=black,xshift=6cm] (0,0) circle (.5ex);
		\draw[black,fill=black,xshift=6cm] (-1,-1) circle (.5ex);
		\draw[red,fill=red,xshift=6cm] (-1,-2) circle (.5ex);
		\draw[fill=black,xshift=6cm] (1,-1) circle (.5ex);
		\draw [thick ,xshift=6cm] (0,0) -- (-1,-1) ;
		\draw [thick,red,xshift=6cm] (-1,-1) -- (-1,-2) ;
		\draw [thick,xshift=6cm] (0,0) -- (1,-1) ;

		\draw[black,fill=black,xshift=9cm] (0,0) circle (.5ex);
		\draw[red,fill=red,xshift=9cm] (-1,-1) circle (.5ex);
		\draw[red,fill=red,xshift=9cm] (1,-2) circle (.5ex);
		\draw[black,fill=black,xshift=9cm] (1,-1) circle (.5ex);
		\draw [thick, red,xshift=9cm] (0,0) -- (-1,-1) ;
		\draw [thick,xshift=9cm] (0,0) -- (1,-1) ;
		\draw [thick,red,xshift=9cm] (1,-1) -- (1,-2) ;
		
		\draw[black,fill=black,xshift=12cm] (0,0) circle (.5ex);
		\draw[black,fill=black,xshift=12cm] (0,-1) circle (.5ex);
		\draw[black,fill=black,xshift=12cm] (1,-2) circle (.5ex);
		\draw[red,fill=red,xshift=12cm] (-1,-2) circle (.5ex);
		\draw [thick, xshift=12cm] (0,0) -- (0,-1) ;
		\draw [thick,xshift=12cm] (0,-1) -- (1,-2) ;
		\draw [thick,red,xshift=12cm] (0,-1) -- (-1,-2) ;

		\draw[black,fill=black,xshift=15cm] (0,0) circle (.5ex);
		\draw[black,fill=black,xshift=15cm] (0,-1) circle (.5ex);
		\draw[red,fill=red,xshift=15cm] (-1,-1) circle (.5ex);
		\draw[fill=black,xshift=15cm] (1,-1) circle (.5ex);
		\draw [thick,red, xshift=15cm] (0,0) -- (-1,-1) ;
		\draw [thick,xshift=15cm] (0,0) -- (1,-1) ;
		\draw [thick,xshift=15cm] (0,0) -- (0,-1) ;

	\end{tikzpicture}
	\caption{All 5 ordered trees with 4 nodes and the old leaves marked in red.}
\end{center}
\end{figure}

Just for interest, we sketch the basic generating function, involving the count of edges, old resp.\ young leaves, using the three variables $z, u, v$ in that order:
\begin{align*}
	F=F(z,u,v)&=1+z(F-1+u)\sum_{k\ge0}(z(F-1+v))^{k}\\*
	&=1+\frac{z(F-1+u)}{1- z(F-1+u)},
\end{align*}
hence
\begin{equation*}
	F=\frac{1+z-uz-\sqrt{1-2(1+v)z+((1+v)^2-4u)z^2}}{2z}.
\end{equation*}
If $v=u$ (no distinction between old resp.\ young), the equation simplifies to 
\begin{equation*}
	F=\frac{1+z-uz-\sqrt{1-2(1+u)z+(1-u)^2z^2}}{2z}.
\end{equation*}

Chen, Deutsch, and Elizalde~\cite{CDE}, using the Lagrange inversion formula, manage to extract the coefficients of
$z^nu^iv^j$ ($n$ edges, $i$ old leaves, $j$ young leaves) in nice form, and derive many further properties.
The instance of $v=u$ leads to the enumeration of \emph{all} leaves and thus to the celebrated \emph{Narayana} numbers.

A related concept is as follows: instead of considering leaves (or equivalently, the edges leading to a leaf), we consider
\emph{all} edges and single out those edges that are the leftmost edges of its parent. A small example explains this:
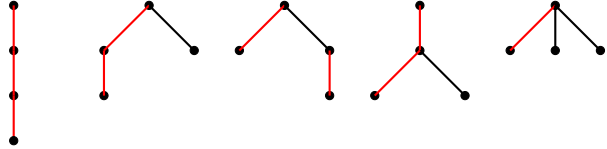
\begin{figure}[h]
		\begin{center}
	\begin{tikzpicture}[scale=0.6]

		\draw[black,fill=black,xshift=3cm] (0,0) circle (.5ex);
		\draw[black,fill=black,xshift=3cm] (0,-1) circle (.5ex);
		\draw[black,fill=black,xshift=3cm] (0,-2) circle (.5ex);
		
		\draw [ thick, red,xshift=3cm] (0,0) -- (0,-1) ;
		\draw [thick,red,xshift=3cm] (0,-1) -- (0,-2) ;
		\draw [thick,red,xshift=3cm] (0,-2) -- (0,-3) ;
		\draw[fill=black,xshift=3cm] (0,-3) circle (.5ex);
		
		\draw[black,fill=black,xshift=6cm] (0,0) circle (.5ex);
		\draw[black,fill=black,xshift=6cm] (-1,-1) circle (.5ex);
		\draw[fill=black,xshift=6cm] (-1,-2) circle (.5ex);
		\draw[fill=black,xshift=6cm] (1,-1) circle (.5ex);
		\draw [thick ,red,xshift=6cm] (0,0) -- (-1,-1) ;
		\draw [thick,red,xshift=6cm] (-1,-1) -- (-1,-2) ;
		\draw [thick,xshift=6cm] (0,0) -- (1,-1) ;

		\draw[black,fill=black,xshift=9cm] (0,0) circle (.5ex);
		\draw[fill=black,xshift=9cm] (-1,-1) circle (.5ex);
		\draw[fill=black,xshift=9cm] (1,-2) circle (.5ex);
		\draw[black,fill=black,xshift=9cm] (1,-1) circle (.5ex);
		\draw [thick, red,xshift=9cm] (0,0) -- (-1,-1) ;
		\draw [thick,xshift=9cm] (0,0) -- (1,-1) ;
		\draw [thick,red,xshift=9cm] (1,-1) -- (1,-2) ;
		
		\draw[black,fill=black,xshift=12cm] (0,0) circle (.5ex);
		\draw[black,fill=black,xshift=12cm] (0,-1) circle (.5ex);
		\draw[black,fill=black,xshift=12cm] (1,-2) circle (.5ex);
		\draw[fill=black,xshift=12cm] (-1,-2) circle (.5ex);
		\draw [thick, red,xshift=12cm] (0,0) -- (0,-1) ;
		\draw [thick,xshift=12cm] (0,-1) -- (1,-2) ;
		\draw [thick,red,xshift=12cm] (0,-1) -- (-1,-2) ;

		\draw[black,fill=black,xshift=15cm] (0,0) circle (.5ex);
		\draw[black,fill=black,xshift=15cm] (0,-1) circle (.5ex);
		\draw[ fill=black,xshift=15cm] (-1,-1) circle (.5ex);
		\draw[fill=black,xshift=15cm] (1,-1) circle (.5ex);
		\draw [thick,red, xshift=15cm] (0,0) -- (-1,-1) ;
		\draw [thick,xshift=15cm] (0,0) -- (1,-1) ;
		\draw [thick,xshift=15cm] (0,0) -- (0,-1) ;

	\end{tikzpicture}
	\caption{All 5 ordered trees with 4 nodes and each leftmost edge marked in red.}
\end{center}
\end{figure}
The enumeration with respect to the number of nodes ($z$), old edges ($u$) and young edges ($v$) is not hard.
The coefficient of $z^nu^iv^j$ is of interest; note that for $n$ nodes, the number of young edges is $j=n-1-i$, so we can ignore the variable $v$.

The generating function $F(z,u)$ is then
\begin{equation*}
	F(z,u)=z+zu\frac{F(z,u)}{1-F(z,u)}=z\Phi(F(z,u)),
\end{equation*}
with $\Phi(y)=1+\frac{uy}{1-y}$, where we consider $u$ to be a dummy parameter. We can use the Lagrange inversion formula:
\begin{align*}
	[z^n]F(z,u)=\frac1n[y^{n-1}]\Phi(y)^n
\end{align*}
and further 
\begin{align*}
	[z^nu^k]F(z,u)&=\frac1n[y^{n-1}][u^k]\Big(1+\frac{uy}{1-y}\Big)^n\\
	&=\frac1n[y^{n-1}]\binom nk\Big(\frac{y}{1-y}\Big)^{k}=\frac1n[y^{n-1-nk}]\binom nk\Big(\frac{1}{1-y}\Big)^{k}\\
	&=\frac1n\binom nk\binom{n-1}{n-1-kn}=\frac1n\binom nk\binom{n-1}{kn}.
\end{align*}
We may note as well that
\begin{equation*}
	F(z,u)=\frac{1+z(1-u)-\sqrt{1-2(1+u)z+(1-u)^2z^2}}{2}.
\end{equation*}

After these introductory and motivating remarks we now come to the main tree model in this paper: instead of considering the (unfair) model of one edge being old and the others being young we
consider a \emph{fair} model: exactly \emph{one} edge emanating from each node that is not a leaf is \emph{distinguished}. The total number of such trees with
$n$ nodes is than naturally higher than the usual Catalan number.
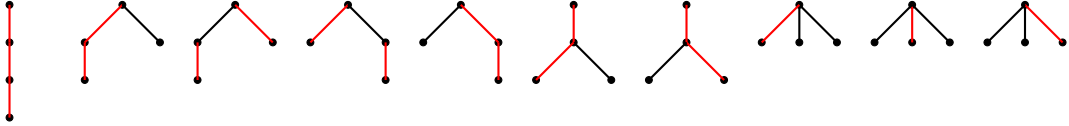
\begin{figure}[h]
	\begin{center}
	\begin{tikzpicture}[scale=0.5]

		\draw[black,fill=black,xshift=3cm] (0,0) circle (.5ex);
		\draw[black,fill=black,xshift=3cm] (0,-1) circle (.5ex);
		\draw[black,fill=black,xshift=3cm] (0,-2) circle (.5ex);
		\draw[black,fill=black,xshift=3cm] (0,-3) circle (.5ex);
		\draw [ thick,red, xshift=3cm] (0,0) -- (0,-1) ;
		\draw [thick,red,xshift=3cm] (0,-1) -- (0,-2) ;
		\draw [thick,red,xshift=3cm] (0,-2) -- (0,-3) ;
		
		\draw[black,fill=black,xshift=6cm] (0,0) circle (.5ex);
		\draw[black,fill=black,xshift=6cm] (-1,-1) circle (.5ex);
		\draw[black
		,fill=black,xshift=6cm] (-1,-2) circle (.5ex);
		\draw[black,fill=black,xshift=6cm] (1,-1) circle (.5ex);
		\draw [thick,red ,xshift=6cm] (0,0) -- (-1,-1) ;
		\draw [thick,red,xshift=6cm] (-1,-1) -- (-1,-2) ;
		\draw [thick,xshift=6cm] (0,0) -- (1,-1) ;
		
		\draw[black,fill=black,xshift=9cm] (0,0) circle (.5ex);
		\draw[black,fill=black,xshift=9cm] (-1,-1) circle (.5ex);
		\draw[black,fill=black,xshift=9cm] (-1,-2) circle (.5ex);
		\draw[black,fill=black,xshift=9cm] (1,-1) circle (.5ex);
		\draw [thick, xshift=9cm] (0,0) -- (-1,-1) ;
		\draw [thick,red,xshift=9cm] (-1,-1) -- (-1,-2) ;
		\draw [thick,red,xshift=9cm] (0,0) -- (1,-1) ;

		\draw[black,fill=black,xshift=12cm] (0,0) circle (.5ex);
		\draw[black,fill=black,xshift=12cm] (-1,-1) circle (.5ex);
		\draw[black,fill=black,xshift=12cm] (1,-2) circle (.5ex);
		\draw[black,fill=black,xshift=12cm] (1,-1) circle (.5ex);
		\draw [thick, red,xshift=12cm] (0,0) -- (-1,-1) ;
		\draw [thick,xshift=12cm] (0,0) -- (1,-1) ;
		\draw [thick,red,xshift=12cm] (1,-1) -- (1,-2) ;
		
		\draw[black,fill=black,xshift=15cm] (0,0) circle (.5ex);
		\draw[black,fill=black,xshift=15cm] (-1,-1) circle (.5ex);
		\draw[black,fill=black,xshift=15cm] (1,-2) circle (.5ex);
		\draw[black,fill=black,xshift=15cm] (1,-1) circle (.5ex);
		\draw [thick, xshift=15cm] (0,0) -- (-1,-1) ;
		\draw [thick,red,xshift=15cm] (0,0) -- (1,-1) ;
		\draw [thick,red,xshift=15cm] (1,-1) -- (1,-2) ;

		\draw[black,fill=black,xshift=18cm] (0,0) circle (.5ex);
		\draw[black,fill=black,xshift=18cm] (0,-1) circle (.5ex);
		\draw[black,fill=black,xshift=18cm] (1,-2) circle (.5ex);
		\draw[black,fill=black,xshift=18cm] (-1,-2) circle (.5ex);
		\draw [thick, red,xshift=18cm] (0,0) -- (0,-1) ;
		\draw [thick,xshift=18cm] (0,-1) -- (1,-2) ;
		\draw [thick,red,xshift=18cm] (0,-1) -- (-1,-2) ;

		\draw[black,fill=black,xshift=21cm] (0,0) circle (.5ex);
		\draw[black,fill=black,xshift=21cm] (0,-1) circle (.5ex);
		\draw[black,fill=black,xshift=21cm] (1,-2) circle (.5ex);
		\draw[black,fill=black,xshift=21cm] (-1,-2) circle (.5ex);
		\draw [thick, red,xshift=21cm] (0,0) -- (0,-1) ;
		\draw [thick,red, xshift=21cm] (0,-1) -- (1,-2) ;
		\draw [thick,xshift=21cm] (0,-1) -- (-1,-2) ;

		\draw[black,fill=black,xshift=24cm] (0,0) circle (.5ex);
		\draw[black,fill=black,xshift=24cm] (0,-1) circle (.5ex);
		\draw[black,fill=black,xshift=24cm] (1,-1) circle (.5ex);
		\draw[black,fill=black,xshift=24cm] (-1,-1) circle (.5ex);
		\draw [thick, xshift=24cm] (0,0) -- (0,-1) ;
		\draw [thick,xshift=24cm] (0,0) -- (1,-1) ;
		\draw [thick,red,xshift=24cm] (0,0) -- (-1,-1) ;
		
		\draw[black,fill=black,xshift=27cm] (0,0) circle (.5ex);
		\draw[black,fill=black,xshift=27cm] (0,-1) circle (.5ex);
		\draw[black,fill=black,xshift=27cm] (1,-1) circle (.5ex);
		\draw[black,fill=black,xshift=27cm] (-1,-1) circle (.5ex);
		\draw [thick,red, xshift=27cm] (0,0) -- (0,-1) ;
		\draw [thick,xshift=27cm] (0,0) -- (1,-1) ;
		\draw [thick,xshift=27cm] (0,0) -- (-1,-1) ;
		
		\draw[black,fill=black,xshift=30cm] (0,0) circle (.5ex);
		\draw[black,fill=black,xshift=30cm] (0,-1) circle (.5ex);
		\draw[black,fill=black,xshift=30cm] (1,-1) circle (.5ex);
		\draw[black,fill=black,xshift=30cm] (-1,-1) circle (.5ex);
		\draw [thick, xshift=30cm] (0,0) -- (0,-1) ;
		\draw [thick,red,xshift=30cm] (0,0) -- (1,-1) ;
		\draw [thick,xshift=30cm] (0,0) -- (-1,-1) ;

	\end{tikzpicture}
	\caption{All 10 ordered trees with 4 nodes and exactly one (distinguished) child marked in red.}
\end{center}
\end{figure}

The rest of this paper is devoted to a thorough analysis of such trees with distinguished children.

\section{Trees with distinguished children}
We use the generating function $A(z)$ for such trees, where $z$ marks the number of nodes. Then
\begin{equation*}
	A=z+zA\sum_{i\ge0}A^i\sum_{j\ge0}A^j=z+zA\sum_{k\ge0}(k+1)A^k=z+\frac{zA}{(1-A)^2},
\end{equation*}
which falls again into the regime of the Lagrange inversion formula;
\begin{equation*}
	A=z\Phi(A),\, \quad \Phi(A)= 1+\frac{A}{(1-A)^2}
\end{equation*}
and therefore
\begin{align*}
	[z^n]A&=\frac1n[A^{n-1}]\Phi^n(A)	
	=\frac1n[A^{n-1}]\Big(1+\frac{A}{(1-A)^2}\Big)^n	\\
	&=\frac1n[A^{n-1}]\sum_{k=0}^{n} \binom{n}{k}\frac{A^k}{(1-A)^{2k}}	
	=\frac1n\sum_{k=0}^{n} \binom{n}{k}[A^{n-1-k}]\frac{1}{(1-A)^{2k}}	\\
	&=\frac1n\sum_{k=0}^{n} \binom{n}{k}\binom{n-2+k}{n-1-k}	
	=\frac1n\sum_{k=1}^{n} \binom{n}{k}\binom{2n-2-k}{k-1}.
\end{align*}
The sequence is $1, 1, 3, 10, 37, 146, 602, 2563, 11181, 49720, 224540$, which is A109081 in \cite{OEIS}; the description contains many references as well.

We are also interested to additionally label leaves with an additional variable $u$; then
\begin{equation*}
	A=z\Phi(A),\, \quad \Phi(A)= u+\frac{A}{(1-A)^2}
\end{equation*}
and
\begin{align*}
	[z^n]A&=\frac1n[A^{n-1}]\Phi^n(A)	
	=\frac1n[A^{n-1}]\Big(u+\frac{A}{(1-A)^2}\Big)^n	\\
	&=\frac1n[A^{n-1}]\sum_{k=0}^{n} \binom{n}{k}u^{n-k}\frac{A^k}{(1-A)^{2k}}	
	=\frac1n\sum_{k=0}^{n} \binom{n}{k}u^{n-k}[A^{n-1-k}]\frac{1}{(1-A)^{2k}}	\\
	&=\frac1n\sum_{k=0}^{n} \binom{n}{k}u^{n-k}\binom{n-2+k}{n-1-k}	
	=\frac1n\sum_{k=1}^{n} u^k\binom{n}{k}\binom{2n-2-k}{k-1};
\end{align*}
further
\begin{align*}
	[z^nu^j]A&=\frac1n[A^{n-1}]\Phi^n(A)	\\
	&=[u^j]\frac1n\sum_{k=1}^{n} u^k\binom{n}{k}\binom{2n-2-k}{k-1}	
	=\frac1n\binom{n}{j}\binom{2n-2-j}{j-1}.
\end{align*}

Now we move to the asymptotics of this sequence. The individual steps can be seen following the analysis of simply generated families of tree, see e.\,g., \cite{DrmotaBook}:
\begin{equation*}
	A=z\Phi(A), \quad \Phi(A)= 1+\frac{A}{(1-A)^2},
\end{equation*}
\begin{equation*}
	\Phi(\tau)=\tau \Phi'(\tau), \quad -1+3\tau-\tau^2+\tau^3=0,
\end{equation*}
\begin{equation*}
	\tau=	0.3611030805;\quad\rho=\frac{\tau}{\Phi(\tau)}=0.1916025857.
\end{equation*}
We can throw in another variable, and consider
\begin{equation*}
	A=zu+\frac{zA}{(1-A)^2}.
\end{equation*}
Then the coefficient of $z^nu^k$ in $A(z,u)$ is the number of trees with distinguished children, and $n$ is the number of nodes, and $k$ is the number of leaves.
This is the equivalent of the Narayana generating function related to the tree family. We start with the functional equation
\begin{equation*}
	A=z\Phi(A),\, \quad \Phi(A)= u+\frac{A}{(1-A)^2}
\end{equation*}
and continue
\begin{align*}
	[z^n]A&=\frac1n[A^{n-1}]\Phi^n(A)	
	=\frac1n[A^{n-1}]\Big(u+\frac{A}{(1-A)^2}\Big)^n	\\
	&=\frac1n[A^{n-1}]\sum_{k=0}^{n} \binom{n}{k}u^{n-k}\frac{A^k}{(1-A)^{2k}}	
	=\frac1n\sum_{k=0}^{n} \binom{n}{k}u^{n-k}[A^{n-1-k}]\frac{1}{(1-A)^{2k}}	\\
	&=\frac1n\sum_{k=0}^{n} \binom{n}{k}u^{n-k}\binom{n-2+k}{n-1-k}	
	=\frac1n\sum_{k=1}^{n} u^k\binom{n}{k}\binom{2n-2-k}{k-1}.
\end{align*}
Hence the number of (new) trees with $n$ nodes and $k$ leaves is given by
\begin{align*}
	[z^nu^k]A(z,u)
	=\frac1n\binom{n}{k}\binom{2n-2-k}{k-1}.
\end{align*}
Now we move to the average number of leaves;
\begin{equation*}
	\frac{\partial}{\partial u}A(z,u)=\frac{zA(z,1)^2-2zA(z	,1)+z}{3A(z,1)^2-2(2+zu)A(z	,1)+1-z+2zu}.
\end{equation*}
\begin{equation*}
	B(z):=	\frac{\partial}{\partial u}A(z,u)\Big|_{u=1}=\frac{zA(z,1)^2-2zA(z	,1)+z}{3A(z,1)^2-2(2+z)A(z	,1)+1+z}.
\end{equation*}
The quantity $A(z	,1)$ can be eliminated according to its equation, and one finds
\begin{equation*}
	-z+z ( z+5 ) B ( z) +( -4+20z+4{z}
	^{2}+3{z}^{3} ) B ( z) ^{2}+
	( -4+20z+4{z}^{2}+3{z}^{3} ) B( z	) ^{3}=0.
\end{equation*}
Around $z=\rho$, as is known for simply generated trees, see \cite{DrmotaBook}
\begin{equation*}
	A\sim \tau \sqrt{\frac{2\Phi(\tau)}{\Phi{''}(\tau)}}\sqrt{1-\frac{z}{\rho}}
\end{equation*}
and 
\begin{equation*}
	[z^n]A\sim \sqrt{\frac{\Phi(\tau)}{2\pi\Phi{''}(\tau)}}\rho^{-n}n^{-3/2}.
\end{equation*}
We replace $A$ by $\tau- 0.3646859471\sqrt{1-\dfrac{z}{\rho}}$ in $B(z)$ to get an approximation around $z=\rho$. First
\begin{equation*}
	B(z)=\frac {z ( 6z+2{z}^{2}-6{A}^{2}+10A+4z{A}^{2}-2zA-A{	z}^{2}-4 ) }{-4+20z+4{z}^{2}+3{z}^{3}};
\end{equation*}
we start with denominator
\begin{align*}
	-4+20z+4{z}^{2}+3{z}^{3}&\sim \bigl(1-\frac z\rho\bigr)(-3\rho{z}^{2}-4\rho z-3{\rho}^{2}z-20\rho-4{\rho}^{2}-3{\rho}^{3})\\
	&\sim -4.189050369\bigl(1-\frac z\rho\bigr);
\end{align*}
then we consider the numerator:
\begin{align*}
	z ( 6z+2{z}^{2}-6{A}^{2}+10A+4z{A}^{2}-2zA-A{	z}^{2}-4 ) 
	&\sim -.4052983049\sqrt{1-\frac z\rho}.
\end{align*}
Combined
\begin{align*}
	B(z)&\sim \frac{-.4052983049\sqrt{1-\frac z\rho}}{-4.189050369\bigl(1-\frac z\rho\bigr)}
	\sim 0.096751833769\bigl(1-\frac z\rho\bigr)^{-1/2}.
\end{align*}
We will use the classical transfer formulae, see \cite{FO}
\begin{equation*}
	c(1-z)^{-\beta}\Longrightarrow \frac{c}{\Gamma(\beta)}n^{\beta-1}, \qquad
	c(1-z/\rho)^{-\beta}\Longrightarrow \frac{c}{\Gamma(\beta)}n^{\beta-1}\rho^{-n}.
\end{equation*}
Applying this leads to
\begin{align*}
	\frac{[z^n]B(z)}{[z^n]A(z)}&\sim \frac{[z^n](0.096751833769\bigl(1-\frac z\rho\bigr)^{-1/2}}
	{[z^n]\Bigl(.3611030805-.3646859471\sqrt{1-\frac z\rho}\Bigr)}\\
	&\sim \frac{[z^n]0.096751833769\bigl(1-\frac z\rho\bigr)^{-1/2}}
	{[z^n]\Bigl(-.3646859471\sqrt{1-\frac z\rho}\Bigr)}\\
	&\sim -.2653017878\frac{[z^n](1-z)^{-1/2}}
	{[z^n]\sqrt{1-z}}\\
	&\sim -.2653017878\frac{\frac{1}{\Gamma(1/2)}n^{-1/2}}
	{\frac1{\Gamma(-1/2)}n^{-3/2}}\\
	&\sim 0.5306035756n.
\end{align*}
For classical ordered trees the corresponding answer is $\sim n/2$, so we have slightly more leaves, on average for the new tree model.

Now we move the (average) height of the new tree family. No new considerations are necessary, since the paper by Flajolet and Odlyzko \cite{FO} deals with simply generated families of trees. The relevant quantities are
\begin{equation*}
	\sqrt{\frac{2}{\phi(\tau)\phi{''}(\tau)}}\phi'(\tau)\cdot \sqrt{\pi n}= 1.009922004\sqrt{\pi n};
\end{equation*}
we might compare this with the average height of ordered trees, which is asymptotic to $\sqrt{\pi n}$. So the average height in the new model is slightly higher.

The next parameter of interest is the \emph{degree of the root.} We go back to the  defining equation of $A(z)$,
\begin{equation*}
	A=z+\frac{zA}{(1-A)^2}=z+z\sum_{k\ge1}kA^k
\end{equation*}
and use an additional parameter $w$, to count the degree of the root:
\begin{equation*}
	F(z,w)=z+z\sum_{k\ge1}kw^kA^k=z+z\frac{wA}{(1-wA)^2}.
\end{equation*}
Eliminating $A$, one can find an algebraic equation for $F$ of order 3, which is a bit long; we will only consider the average of our parameter, in other words, we differentiate w.r.t. $w$, and set $w=1$. We obtain:
\begin{equation*}
	z\sum_{k\ge1}k^2A^k=\frac{zA(A+1)}{(1-A)^3}\sim .361103078-1.079796775\sqrt{1-\frac z\rho}.
\end{equation*}
Comparing this with
\begin{equation*}
	A\sim.3611030805-.3646859471\sqrt{1-\frac z\rho}
\end{equation*}
we find that the average degree of the root tends to the constant
\begin{equation*}
	\frac{-1.079796775}{-.3646859471}=2.960894939.
\end{equation*}
In the classical case of ordered trees, that constant is just 2 (see \cite{Kemp}), so it is slightly higher here.

The next parameter to be studied is the \emph{height of the left-most leaf.} For that, we introduce again a second variable $w$, measuring just the height of the left-most leaf;
\begin{equation*}
	F=zw+\sum_{j\ge1}jwzFA^{j-1}=zw+zwF\sum_{j\ge1}jA^{j-1}=zw+zwF\frac1{(1-A)^2};
\end{equation*}
solving
\begin{equation*}
	F=\frac{zw(1-A)^2}{1-zw-2A+A^2}.
\end{equation*}
Differentiating w.r.t. $w$, followed by $w=1$:
\begin{equation*}
	\frac{z(1-A)^4}{1-z-2A+A^2}=\frac{A^2}{z}\sim.6805515397-1.374607951\sqrt{1-\frac z\rho}.
\end{equation*}
Thus  the height of the left-most leaf tends to the constant
\begin{equation*}
	\frac{-1.374607951}{-.3646859471}=3.769292351.
\end{equation*}

Now we move to the \emph{pathlength}, which is the sum of the distances, summed over all nodes in the new family. Again, an additional variable is used to keep track.
\begin{align*}
	F(z,w)&=\sum_{n,k\ge1}p_{n,k}z^nw^k\\
	&=zw+\sum_{j\ge1}j zw\Bigl(\sum_{n,k\ge1}p_{n,k}(zw)^nw^k\Bigr)^j\\
	&=zw+zw\frac{F(zw,w)}{(1-F(zw,w))^2}.
\end{align*}
Compare this with the classical case in \cite{FS}, where there is no square in the denominator. One can use the Lagrange inversion formula to pull out coefficients, since
\begin{equation*}
	F=z\Phi(F),\qquad \Phi(y)=w+w\frac{y}{(1-y)^2},
\end{equation*}
and
\begin{align*}
	[z^n]F=\frac1n[F^{n-1}]\Big(w+w\frac{F}{(1-F)^2}\Bigr)^n=\frac{w^n}n[F^{n-1}]\Big(1+\frac{F}{(1-F)^2}\Bigr)^n;
\end{align*}
something could be written out, but since it is not too pleasant, we won't do it. 

We can compute the average pathlength.
In order to work with implicitly defined functions, we must count nodes; the distance we measure in terms of the number of edges. Then
\begin{align*}
	F(z,u)&=z+\sum_{k\ge1}k z^{k}F^k(zu,u)=z+\frac{zF(zu,u)}{(1-zF(zu,u))^2}.
\end{align*}
Differentiate this w.r.t. $u$, followed by $u:=1$, using 
\begin{equation*}
	\frac{\partial F(z,u)}{\partial u}=Y,\quad \frac{\partial F(xu,u)}{\partial x}\bigg|_{x=z,u=1}=A'(z),
	\quad \frac{\partial F(zu,v)}{\partial v}\bigg|_{v=1,u=1}=Y.
\end{equation*}
Solving and simplifying, we get
\begin{equation*}
	Y=-{\frac {2zA(z)^{2}-2A(z)^{2}+zA(z)+{z}^{2}A(z)+2A(z)+{z}^{2}-2z}
		{3{			z}^{3}+4{z}^{2}+20z-4}};
\end{equation*}
the series expansion start with
\begin{equation*}
	Y=z^2+7z^3+41z^4+230z^5+1261z^6+6824z^7+36627z^8+195504z^9.
\end{equation*}
The denominator factors:
\begin{align*}
	3{			z}^{3}+4{z}^{2}+20z-4&=(z-\rho)(3z^2+(4+3\rho)z+20+4\rho+3\rho^2)\\
	&\sim -\bigl(1-\frac z\rho\bigr)4.189050371.
\end{align*}
The coefficient of $z^n$ is asymptotic to $4.189050371\rho^{-n}$. For the numerator, we find the approximation
$-.3857113818$ thus 
\begin{equation*}
	Y\sim \frac{-.3857113818}{-\bigl(1-\frac z\rho\bigr)4.189050371}=\frac{0.09207609067}{1-\frac z\rho};
\end{equation*}
the coefficient of $z^n$ in $Y$ is asymptotic to $0.09207609067\rho^{-n}$. This has to divided by the total number of objects of size $n$. Recall that this was
\begin{equation*}
	\sim-1.079796775[z^n]	\sqrt{1-\frac z\rho} \sim -1.079796775\frac1{\Gamma({-\frac12})}n^{-3/2}\rho^{-n}=0.3046050464n^{-3/2}\rho^{-n}.
\end{equation*}
Taking quotients, we find the average pathlength to be asymptotic to 
\begin{equation*}
	\frac{ 0.09207609067\rho^{-n}}{0.3046050464\,n^{-3/2}\rho^{-n}}=0.3022802536\,n^{3/2} .
\end{equation*}

Counting of \emph{old leaves} (leftmost children), as in \cite{CDE} for the classical case: In our model of distinguished edges, we look at old leaves (the total number of leaves)  was already considered.   The one node tree has by definition one node, one old leaf and zero young leaves. Recall that
\begin{equation*}
	A(z,u)= zu+zA(z,u)\sum_{k\ge0}(k+1)A(z,u)^k
\end{equation*}
when $u$ marks all leaves; we use two variables $u$, $v$ for old versus young leaves. Then
\begin{equation*}
	A(z,u,v)= zu+z\sum_{k\ge1}kA(z,u,v)(A(z,u,v)-zu+zv)^k=zu+\frac{zA(z,u,v)}{(1-A(z,u,v)+zu-zv)^2}.
\end{equation*}
We leave further analysis of the generating function to the interested reader. It would be more appealing to set $v=1$, so one would only count old leaves and do not care about
the young leaves. Then one would have to analyze
\begin{equation*}
	A(z,u)= zu+\frac{zA(z,u)}{(1-A(z,u)+zu-z)^2}.
\end{equation*}
The Lagrange inversion formula is perhaps not immediately available, but $A(z,u)$ satisfies an algebraic equation of order 3, and one can work out asymptotics for the average number of old leaves. We compute
\begin{align*}
	&\frac{\partial A(z,u)}{\partial u}\Big|_{u=1}\\&=-\frac {z(     A^2  ( z  )   +2z
		A^2  ( z  )    +2   A^{2}  ( z		 )    {z}^{2}-2A  ( z  ) -12A  ( z ) z-5A  ( z  ) {z}^{2}-2{z}^{3}A  ( z		 ) -14z-2{z}^{3}+4  ) }{-4+20z+4{z}^{2}+3{z}^{3}}
\end{align*}
where $A(z)$ stands for $A(z,1)$. As was already computed,
\begin{equation*}
	-4+20z+4{z}^{2}+3{z}^{3}\sim -21.86322464\bigl(1-\frac z\rho\bigr).
\end{equation*}
Together with the expansion of the denominator, 
\begin{equation*}
	\sim\frac{-.1902455753\sqrt{1-\frac z\rho}}{-21.86322464\bigl(1-\frac z\rho\bigr)}=\frac{0.008701624689}{\sqrt{1-\frac z\rho}}.
\end{equation*}
By transfer, the coefficient of $z^n$ is asymptotic to
\begin{equation*}
	0.008701624689\frac1{\Gamma(\frac12)}n^{-1/2}\rho^{-n}.
\end{equation*}
This has to be divided by the asymptotic equivalent of the total number of object of size $n$, viz. $-1.079796775\frac1{\Gamma({-\frac12})}n^{-3/2}\rho^{-n}$, so we get
\begin{equation*}
	\frac{0.008701624689\frac1{\Gamma(\frac12)}n^{-1/2}\rho^{-n}}{-1.079796775\frac1{\Gamma({-\frac12})}n^{-3/2}\rho^{-n}}=
	\frac{0.008701624689\frac1{\Gamma(\frac12)}n}{-1.079796775\frac1{\Gamma({-\frac12})}}=0.01611715258n.
\end{equation*}
This is the asymptotic equivalent of old leaves in the new tree model, assuming that all objects of $n$ nodes are equally likely.

\section{Marked ordered trees}

In \cite{DMS} we find the following variation of ordered trees. Each rightmost edge might be marked or not, if it does not lead to an endnode (leaf).
We depict a marked edge by the red color and draw all of them of size (4 nodes):
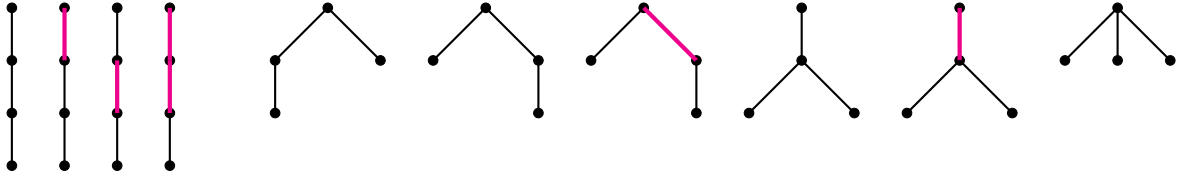
\begin{figure}[h]
	\begin{center}
	\begin{tikzpicture}[scale=0.7]

		\draw[black,fill=black] (0,0) circle (.5ex);
		\draw[black,fill=black] (0,-1) circle (.5ex);
		\draw[black,fill=black] (0,-2) circle (.5ex);
		\draw[black,fill=black] (0,-3) circle (.5ex);
		\draw [thick] (0,0) -- (0,-1)-- (0,-3) ;%
		
		\draw[black,fill=black,xshift=1cm] (0,0) circle (.5ex);
		\draw[black,fill=black,xshift=1cm] (0,-1) circle (.5ex);
		\draw[black,fill=black,xshift=1cm] (0,-2) circle (.5ex);
		\draw[black,fill=black,xshift=1cm] (0,-3) circle (.5ex);
		\draw [ultra thick,magenta,xshift=1cm] (0,0) -- (0,-1) ;
		\draw [thick,xshift=1cm] (0,-1) -- (0,-2) ;
		\draw [thick,xshift=1cm] (0,-2) -- (0,-3) ;
		
		\draw[black,fill=black,xshift=2cm] (0,0) circle (.5ex);
		\draw[black,fill=black,xshift=2cm] (0,-1) circle (.5ex);
		\draw[black,fill=black,xshift=2cm] (0,-2) circle (.5ex);
		\draw[black,fill=black,xshift=2cm] (0,-3) circle (.5ex);
		\draw [thick,xshift=2cm] (0,0) -- (0,-1) ;
		\draw [ultra thick,magenta,xshift=2cm] (0,-1) -- (0,-2) ;
		\draw [thick,xshift=2cm] (0,-2) -- (0,-3) ;
		
		\draw[black,fill=black,xshift=3cm] (0,0) circle (.5ex);
		\draw[black,fill=black,xshift=3cm] (0,-1) circle (.5ex);
		\draw[black,fill=black,xshift=3cm] (0,-2) circle (.5ex);
		\draw[black,fill=black,xshift=3cm] (0,-3) circle (.5ex);
		\draw [ultra thick,magenta, xshift=3cm] (0,0) -- (0,-1) ;
		\draw [ultra thick,magenta,xshift=3cm] (0,-1) -- (0,-2) ;
		\draw [thick,xshift=3cm] (0,-2) -- (0,-3) ;
		
		\draw[black,fill=black,xshift=6cm] (0,0) circle (.5ex);
		\draw[black,fill=black,xshift=6cm] (-1,-1) circle (.5ex);
		\draw[black,fill=black,xshift=6cm] (-1,-2) circle (.5ex);
		\draw[black,fill=black,xshift=6cm] (1,-1) circle (.5ex);
		\draw [thick, xshift=6cm] (0,0) -- (-1,-1) ;
		\draw [thick,xshift=6cm] (-1,-1) -- (-1,-2) ;
		\draw [thick,xshift=6cm] (0,0) -- (1,-1) ;
		
		\draw[black,fill=black,xshift=9cm] (0,0) circle (.5ex);
		\draw[black,fill=black,xshift=9cm] (-1,-1) circle (.5ex);
		\draw[black,fill=black,xshift=9cm] (1,-2) circle (.5ex);
		\draw[black,fill=black,xshift=9cm] (1,-1) circle (.5ex);
		\draw [thick, xshift=9cm] (0,0) -- (-1,-1) ;
		\draw [thick,xshift=9cm] (0,0) -- (1,-1) ;
		\draw [thick,xshift=9cm] (1,-1) -- (1,-2) ;
		
		\draw[black,fill=black,xshift=12cm] (0,0) circle (.5ex);
		\draw[black,fill=black,xshift=12cm] (-1,-1) circle (.5ex);
		\draw[black,fill=black,xshift=12cm] (1,-2) circle (.5ex);
		\draw[black,fill=black,xshift=12cm] (1,-1) circle (.5ex);
		\draw [thick, xshift=12cm] (0,0) -- (-1,-1) ;
		\draw [ultra thick,magenta,xshift=12cm] (0,0) -- (1,-1) ;
		\draw [thick,xshift=12cm] (1,-1) -- (1,-2) ;
		
		\draw[black,fill=black,xshift=15cm] (0,0) circle (.5ex);
		\draw[black,fill=black,xshift=15cm] (0,-1) circle (.5ex);
		\draw[black,fill=black,xshift=15cm] (1,-2) circle (.5ex);
		\draw[black,fill=black,xshift=15cm] (-1,-2) circle (.5ex);
		\draw [thick, xshift=15cm] (0,0) -- (0,-1) ;
		\draw [thick,xshift=15cm] (0,-1) -- (1,-2) ;
		\draw [thick,xshift=15cm] (0,-1) -- (-1,-2) ;
		
		\draw[black,fill=black,xshift=18cm] (0,0) circle (.5ex);
		\draw[black,fill=black,xshift=18cm] (0,-1) circle (.5ex);
		\draw[black,fill=black,xshift=18cm] (1,-2) circle (.5ex);
		\draw[black,fill=black,xshift=18cm] (-1,-2) circle (.5ex);
		\draw [ultra thick, magenta, xshift=18cm] (0,0) -- (0,-1) ;
		\draw [ thick,xshift=18cm] (0,-1) -- (1,-2) ;
		\draw [thick,xshift=18cm] (0,-1) -- (-1,-2) ;
		
		\draw[black,fill=black,xshift=21cm] (0,0) circle (.5ex);
		\draw[black,fill=black,xshift=21cm] (0,-1) circle (.5ex);
		\draw[black,fill=black,xshift=21cm] (1,-1) circle (.5ex);
		\draw[black,fill=black,xshift=21cm] (-1,-1) circle (.5ex);
		\draw [thick, xshift=21cm] (0,0) -- (0,-1) ;
		\draw [thick,xshift=21cm] (0,0) -- (1,-1) ;
		\draw [thick,xshift=21cm] (0,0) -- (-1,-1) ;
		
	\end{tikzpicture}
	\caption{All 10 marked ordered trees with 4 nodes.}
\end{center}
	\label{bull1}
\end{figure}

Now we move to a symbolic equation for the marked ordered trees:

\begin{figure}[h]
		\begin{center}
	\begin{tikzpicture}[scale=1.0,
		s1/.style={circle=10pt,draw=black!90,thick},
		s2/.style={rectangle,draw=black!50,thick},scale=0.5]
		
		\node at ( 7.0,0) { $\mathscr{A}$};
		
		\node at (8.2,0) { $=$};
		\node(c) at (9.7,0)[s1]{};
		\node at (11.5,0) {$+$};

		\node[xshift=5cm](ch) at (5,1)[s1]{};
		\node[xshift=5cm](e) at (3,-1){ $\mathscr{A}$};
		\node[xshift=5cm](ee) at (5,-1){$\cdots$};
		\node[xshift=5cm](f) at (7,-1){ $\mathscr{A}$};
		\path [draw,-,black!90] (ch) -- (e) node{};%
		\path [draw,-,black!90] (ch) -- (f) node{};%
		\node[xshift=5cm](g) at (9,-1){ $\mathscr{A}$};
		\path [draw,-,black,ultra thick] (ch) -- (8.8+10,-0.5) node {};%
		
		\node[xshift=10cm] at (0.7,0) {$+$};

		\node[xshift=10cm](ch) at (5,1)[s1]{};
		\node[xshift=10cm](e) at (3,-1){ $\mathscr{A}$};
		\node[xshift=10cm](ee) at (5,-1){$\cdots$};
		\node[xshift=10cm](f) at (7,-1){ $\mathscr{A}$};
		\path [draw,-,black!90] (ch) -- (e) node{};%
		\path [draw,-,black!90] (ch) -- (f) node{};%
		\node[xshift=10cm](g) at (9,-1){ $\mathscr{A}\setminus\{\circ\}$};
		\path [draw,-,black,magenta,ultra thick] (ch) -- (8.8+20,-0.5) node {};%

	\end{tikzpicture}
	\caption{Symbolic equation for marked ordered trees, $\mathscr{A}\cdots\mathscr{A}$ refers to $\ge0$ copies of $\mathscr{A}$.}
\end{center}
\end{figure}

In terms of generating functions,
\begin{equation*}
	A(z)=z+\frac{z}{1-A(z)}z+\frac{z}{1-A(z)}2(A(z)-z),
\end{equation*}
with the solution
\begin{equation*}
	A(z)=\frac{1-z-\sqrt{1-6z+5z^2}}{2}=z+z^2+z^3+z^3+10z^4+36z^5+\cdots.
\end{equation*}

The importance of this family of trees lies in the bijection to skew Dyck paths, as given in \cite{DMS}. One walks around the tree as one usually does and translates it into a Dyck path.
The only difference are the red edges. On the way down, nothing special is to be reported, but on the way up, it is translated into a skew step $(-1,-1)$. The present author believes that trees are more manageable when it comes to enumeration issues than skew Dyck paths.

The 10 trees of Figure~\ref{bull1} translate as follows:
\begin{equation*}
	\begin{tikzpicture}[scale=0.3]
		\draw (0,0)--(3,3)--(6,0);
	\end{tikzpicture}
	\quad
	\begin{tikzpicture}[scale=0.3]
		\draw (0,0)--(3,3)--(5,1)--(4,0);
	\end{tikzpicture}
	\quad
	\begin{tikzpicture}[scale=0.3]
		\draw (0,0)--(3,3)--(4,2)--(3,1)--(4,0);
	\end{tikzpicture}
	\quad
	\begin{tikzpicture}[scale=0.3]
		\draw (0,0)--(3,3)--(4,2)--(3,1)--(2,0);
	\end{tikzpicture}
	\quad
	\begin{tikzpicture}[scale=0.3]
		\draw (0,0)--(2,2)--(4,0)--(5,1)--(6,0);
	\end{tikzpicture}
\end{equation*}
\begin{equation*}
	\begin{tikzpicture}[scale=0.3]
		\draw (0,0)--(1,1)--(2,0)--(4,2)--(6,0);
	\end{tikzpicture}
	\quad
	\begin{tikzpicture}[scale=0.3]
		\draw (0,0)--(1,1)--(2,0)--(4,2)--(5,1)--(4,0);
	\end{tikzpicture}
	\quad
	\begin{tikzpicture}[scale=0.3]
		\draw (0,0)--(1,1)--(2,2)--(3,1)--(4,2)--(6,0);
	\end{tikzpicture}
	\quad
	\begin{tikzpicture}[scale=0.3]
		\draw (0,0)--(1,1)--(2,2)--(3,1)--(4,2)--(5,1)--(4,0);
	\end{tikzpicture}
	\quad
	\begin{tikzpicture}[scale=0.3]
		\draw (0,0)--(1,1)--(2,0)--(3,1)--(4,0)--(5,1)--(6,0);
	\end{tikzpicture}
\end{equation*}

Now we combine the concepts `marked ordered trees' and `distinguished children' by distinguish exactly one child of each node in a marked ordered tree.
First, to clarify, we provide an example.
\begin{figure}[h]
		\begin{center}
	\begin{tikzpicture}[scale=0.8]

		\draw[black,fill=black] (0,0) circle (.5ex);
		\draw[black,fill=black] (0,-1) circle (.5ex);
		\draw[black,fill=black] (0,-2) circle (.5ex);
		\draw[black,fill=black] (0,-3) circle (.5ex);
		\draw [ultra thick, densely dotted] (0,0) -- (0,-1)-- (0,-3) ;%
		
		\draw[black,fill=black,xshift=1cm] (0,0) circle (.5ex);
		\draw[black,fill=black,xshift=1cm] (0,-1) circle (.5ex);
		\draw[black,fill=black,xshift=1cm] (0,-2) circle (.5ex);
		\draw[black,fill=black,xshift=1cm] (0,-3) circle (.5ex);
		\draw [ultra thick,magenta,ultra thick, densely dotted, xshift=1cm] (0,0) -- (0,-1) ;
		\draw [thick,ultra thick, densely dotted, xshift=1cm] (0,-1) -- (0,-2) ;
		\draw [thick,ultra thick, densely dotted, xshift=1cm] (0,-2) -- (0,-3) ;
		
		\draw[black,fill=black,xshift=2cm] (0,0) circle (.5ex);
		\draw[black,fill=black,xshift=2cm] (0,-1) circle (.5ex);
		\draw[black,fill=black,xshift=2cm] (0,-2) circle (.5ex);
		\draw[black,fill=black,xshift=2cm] (0,-3) circle (.5ex);
		\draw [thick,ultra thick, densely dotted,xshift=2cm] (0,0) -- (0,-1) ;
		\draw [ultra thick,ultra thick, densely dotted,magenta,xshift=2cm] (0,-1) -- (0,-2) ;
		\draw [thick,ultra thick, densely dotted,xshift=2cm] (0,-2) -- (0,-3) ;
		
		\draw[black,fill=black,xshift=3cm] (0,0) circle (.5ex);
		\draw[black,fill=black,xshift=3cm] (0,-1) circle (.5ex);
		\draw[black,fill=black,xshift=3cm] (0,-2) circle (.5ex);
		\draw[black,fill=black,xshift=3cm] (0,-3) circle (.5ex);
		\draw [ultra thick,ultra thick, densely dotted, magenta, xshift=3cm] (0,0) -- (0,-1) ;
		\draw [ultra thick,ultra thick, densely dotted, magenta,xshift=3cm] (0,-1) -- (0,-2) ;
		\draw [thick,ultra thick, densely dotted, xshift=3cm] (0,-2) -- (0,-3) ;
		
		\draw[black,fill=black,xshift=6cm] (0,0) circle (.5ex);
		\draw[black,fill=black,xshift=6cm] (-1,-1) circle (.5ex);
		\draw[black,fill=black,xshift=6cm] (-1,-2) circle (.5ex);
		\draw[black,fill=black,xshift=6cm] (1,-1) circle (.5ex);
		\draw [thick,densely dotted, xshift=6cm] (0,0) -- (-1,-1) ;
		\draw [thick,densely dotted,xshift=6cm] (-1,-1) -- (-1,-2) ;
		\draw [thick,xshift=6cm] (0,0) -- (1,-1) ;
		
		\draw[black,fill=black,xshift=6cm] (0,0) circle (.5ex);
		\draw[black,fill=black,xshift=6cm] (-1,-1) circle (.5ex);
		\draw[black,fill=black,xshift=6cm] (-1,-2) circle (.5ex);
		\draw[black,fill=black,xshift=6cm] (1,-1) circle (.5ex);
		\draw [thick, ultra thick, densely dotted,xshift=6cm] (0,0) -- (-1,-1) ;
		\draw [thick,ultra thick, densely dotted,xshift=6cm] (-1,-1) -- (-1,-2) ;
		\draw [thick,xshift=6cm] (0,0) -- (1,-1) ;
		
		\draw[black,fill=black,xshift=9cm] (0,0) circle (.5ex);
		\draw[black,fill=black,xshift=9cm] (-1,-1) circle (.5ex);
		\draw[black,fill=black,xshift=9cm] (-1,-2) circle (.5ex);
		\draw[black,fill=black,xshift=9cm] (1,-1) circle (.5ex);
		\draw [thick, ultra thick,xshift=9cm] (0,0) -- (-1,-1) ;
		\draw [thick,ultra thick, densely dotted,xshift=9cm] (-1,-1) -- (-1,-2) ;
		\draw [ultra thick, densely dotted,xshift=9cm] (0,0) -- (1,-1) ;

		\draw[black,fill=black,xshift=12cm] (0,0) circle (.5ex);
		\draw[black,fill=black,xshift=12cm] (-1,-1) circle (.5ex);
		\draw[black,fill=black,xshift=12cm] (1,-2) circle (.5ex);
		\draw[black,fill=black,xshift=12cm] (1,-1) circle (.5ex);
		\draw [thick, ultra thick, densely dotted,xshift=12cm] (0,0) -- (-1,-1) ;
		\draw [thick,xshift=12cm] (0,0) -- (1,-1) ;
		\draw [thick,ultra thick, densely dotted,xshift=12cm] (1,-1) -- (1,-2) ;
		
		\draw[black,fill=black,xshift=15cm] (0,0) circle (.5ex);
		\draw[black,fill=black,xshift=15cm] (-1,-1) circle (.5ex);
		\draw[black,fill=black,xshift=15cm] (1,-2) circle (.5ex);
		\draw[black,fill=black,xshift=15cm] (1,-1) circle (.5ex);
		\draw [thick, xshift=15cm] (0,0) -- (-1,-1) ;
		\draw [thick,ultra thick, densely dotted,xshift=15cm] (0,0) -- (1,-1) ;
		\draw [thick,ultra thick, densely dotted,xshift=15cm] (1,-1) -- (1,-2) ;

		\draw[black,fill=black,xshift=0cm,yshift=-4cm] (0,0) circle (.5ex);
		\draw[black,fill=black,xshift=0cm,yshift=-4cm] (-1,-1) circle (.5ex);
		\draw[black,fill=black,xshift=0cm,yshift=-4cm] (1,-2) circle (.5ex);
		\draw[black,fill=black,xshift=0cm,yshift=-4cm] (1,-1) circle (.5ex);
		\draw [thick,ultra thick, densely dotted, xshift=0cm,yshift=-4cm] (0,0) -- (-1,-1) ;
		\draw [ultra thick,magenta,xshift=0cm,yshift=-4cm] (0,0) -- (1,-1) ;
		\draw [thick,ultra thick, densely dotted,xshift=0cm,yshift=-4cm] (1,-1) -- (1,-2) ;
		
		\draw[black,fill=black,xshift=3cm,yshift=-4cm] (0,0) circle (.5ex);
		\draw[black,fill=black,xshift=3cm,yshift=-4cm] (-1,-1) circle (.5ex);
		\draw[black,fill=black,xshift=3cm,yshift=-4cm] (1,-2) circle (.5ex);
		\draw[black,fill=black,xshift=3cm,yshift=-4cm] (1,-1) circle (.5ex);
		\draw [thick, xshift=3cm,yshift=-4cm] (0,0) -- (-1,-1) ;
		\draw [ultra thick,ultra thick, densely dotted,magenta,xshift=3cm,yshift=-4cm] (0,0) -- (1,-1) ;
		\draw [thick,ultra thick, densely dotted,xshift=3cm,yshift=-4cm] (1,-1) -- (1,-2) ;
		
		\draw[black,fill=black,xshift=6cm,yshift=-4cm] (0,0) circle (.5ex);
		\draw[black,fill=black,xshift=6cm,yshift=-4cm] (0,-1) circle (.5ex);
		\draw[black,fill=black,xshift=6cm,yshift=-4cm] (1,-2) circle (.5ex);
		\draw[black,fill=black,xshift=6cm,yshift=-4cm] (-1,-2) circle (.5ex);
		\draw [thick, ultra thick, densely dotted,xshift=6cm,yshift=-4cm] (0,0) -- (0,-1) ;
		\draw [thick,xshift=6cm,yshift=-4cm] (0,-1) -- (1,-2) ;
		\draw [thick,ultra thick, densely dotted,xshift=6cm,yshift=-4cm] (0,-1) -- (-1,-2) ;
		
		\draw[black,fill=black,xshift=9cm,yshift=-4cm] (0,0) circle (.5ex);
		\draw[black,fill=black,xshift=9cm,yshift=-4cm] (0,-1) circle (.5ex);
		\draw[black,fill=black,xshift=9cm,yshift=-4cm] (1,-2) circle (.5ex);
		\draw[black,fill=black,xshift=9cm,yshift=-4cm] (-1,-2) circle (.5ex);
		\draw [thick, ultra thick, densely dotted,xshift=9cm,yshift=-4cm] (0,0) -- (0,-1) ;
		\draw [thick,ultra thick, densely dotted,xshift=9cm,yshift=-4cm] (0,-1) -- (1,-2) ;
		\draw [thick,xshift=9cm,yshift=-4cm] (0,-1) -- (-1,-2) ;
		
		\draw[black,fill=black,xshift=12cm,yshift=-4cm] (0,0) circle (.5ex);
		\draw[black,fill=black,xshift=12cm,yshift=-4cm] (0,-1) circle (.5ex);
		\draw[black,fill=black,xshift=12cm,yshift=-4cm] (1,-2) circle (.5ex);
		\draw[black,fill=black,xshift=12cm,yshift=-4cm] (-1,-2) circle (.5ex);
		\draw [ultra thick, ultra thick, densely dotted,magenta, xshift=12cm,yshift=-4cm] (0,0) -- (0,-1) ;
		\draw [ thick,xshift=12cm,yshift=-4cm] (0,-1) -- (1,-2) ;
		\draw [thick,ultra thick, densely dotted,xshift=12cm,yshift=-4cm] (0,-1) -- (-1,-2) ;
		
		\draw[black,fill=black,xshift=15cm,yshift=-4cm] (0,0) circle (.5ex);
		\draw[black,fill=black,xshift=15cm,yshift=-4cm] (0,-1) circle (.5ex);
		\draw[black,fill=black,xshift=15cm,yshift=-4cm] (1,-2) circle (.5ex);
		\draw[black,fill=black,xshift=15cm,yshift=-4cm] (-1,-2) circle (.5ex);
		\draw [ultra thick, ultra thick, densely dotted,magenta, xshift=15cm,yshift=-4cm] (0,0) -- (0,-1) ;
		\draw [ thick,ultra thick, densely dotted,xshift=15cm,yshift=-4cm] (0,-1) -- (1,-2) ;
		\draw [thick,xshift=15cm,yshift=-4cm] (0,-1) -- (-1,-2) ;

		\draw[black,fill=black,xshift=0cm,yshift=-7cm] (0,0) circle (.5ex);
		\draw[black,fill=black,xshift=0cm,yshift=-7cm] (0,-1) circle (.5ex);
		\draw[black,fill=black,xshift=0cm,yshift=-7cm] (1,-1) circle (.5ex);
		\draw[black,fill=black,xshift=0cm,yshift=-7cm] (-1,-1) circle (.5ex);
		\draw [thick, xshift=0cm,yshift=-7cm] (0,0) -- (0,-1) ;
		\draw [thick,xshift=0cm,yshift=-7cm] (0,0) -- (1,-1) ;
		\draw [thick,ultra thick, densely dotted,xshift=0cm,yshift=-7cm] (0,0) -- (-1,-1) ;
		
		\draw[black,fill=black,xshift=3cm,yshift=-7cm] (0,0) circle (.5ex);
		\draw[black,fill=black,xshift=3cm,yshift=-7cm] (0,-1) circle (.5ex);
		\draw[black,fill=black,xshift=3cm,yshift=-7cm] (1,-1) circle (.5ex);
		\draw[black,fill=black,xshift=3cm,yshift=-7cm] (-1,-1) circle (.5ex);
		\draw [thick,ultra thick, densely dotted, xshift=3cm,yshift=-7cm] (0,0) -- (0,-1) ;
		\draw [thick,xshift=3cm,yshift=-7cm] (0,0) -- (1,-1) ;
		\draw [thick,xshift=3cm,yshift=-7cm] (0,0) -- (-1,-1) ;
		
		\draw[black,fill=black,xshift=6cm,yshift=-7cm] (0,0) circle (.5ex);
		\draw[black,fill=black,xshift=6cm,yshift=-7cm] (0,-1) circle (.5ex);
		\draw[black,fill=black,xshift=6cm,yshift=-7cm] (1,-1) circle (.5ex);
		\draw[black,fill=black,xshift=6cm,yshift=-7cm] (-1,-1) circle (.5ex);
		\draw [thick, xshift=6cm,yshift=-7cm] (0,0) -- (0,-1) ;
		\draw [thick,ultra thick, densely dotted,xshift=6cm,yshift=-7cm] (0,0) -- (1,-1) ;
		\draw [thick,xshift=6cm,yshift=-7cm] (0,0) -- (-1,-1) ;
		
	\end{tikzpicture}
	\caption{All 17 marked ordered trees with 4 nodes and one edge distinguished (dotted)}
\end{center}
\end{figure}
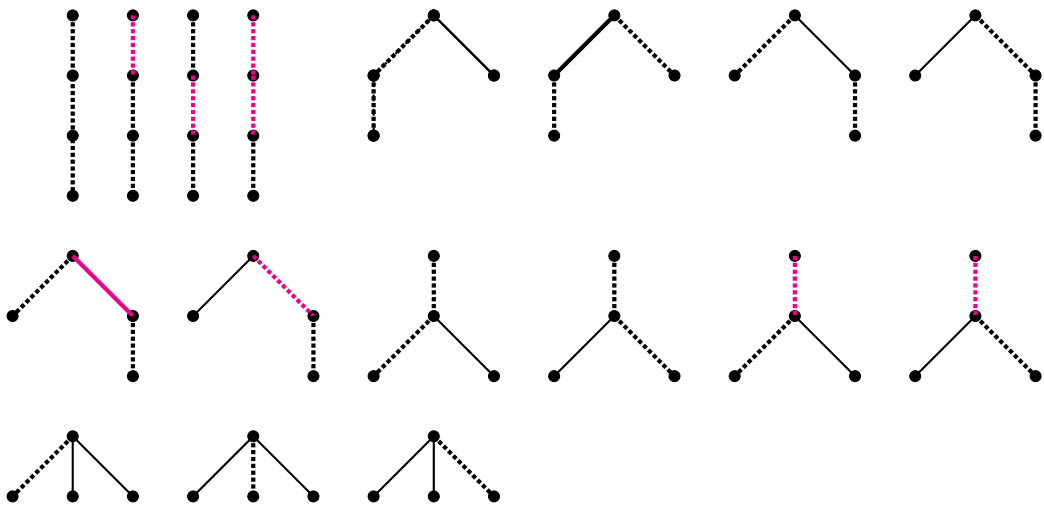

As before, we start from 
\begin{equation*}
	A(z)=z+\sum_{k\ge1}A(z)^k+\sum_{k\ge1}A(z)^{k-1}(A(z)-z)
\end{equation*}
from the old model and get 
\begin{equation*}
	A(z)=z+\sum_{k\ge1}\mathbf{k}A(z)^k+\sum_{k\ge1}\mathbf{k}A(z)^{k-1}(A(z)-z)
\end{equation*}
in the new model. This means
\begin{equation*}
	A=z+\frac{zA}{(1-A)^2}+\frac{z(A-z)}{(1-A)^2}.
\end{equation*}
The algebraic equation of interest is $-A^3+2A^2-A+zA^2+z-z^2=0$; the power series expansion starts as
\begin{equation*}
	A(z)=z+z^2+4 z^3+\mathbf{17} z^4+78 z^5+378 z^6+1906 z^7+9901 z^8+52630 z^9+284926 z^{10}+\cdots
\end{equation*}
The number 17 (in bold) refers to the 17 objects of 4 vertices in the previous example.

Further investigations of parameters in this interesting tree model (marked, distinguished) are left to interested reader and/or future research. The sequence $1,1,4,17,78,\dots$ is not in the encyclopedia of integer sequences.

\section{Concluding remarks}

The major contribution of this paper is the introduction of the model of distinguished children and its variants. The author hopes that interested readers will pick up the subject and contribute some more elaborate analysis on various aspects.

\vspace{24pt}
\noindent {\bf\Large Acknowledgments}\vspace*{6pt}\\
Friendly comments of a reviewer are gratefully acknowledged. This research was not funded.

\bibliographystyle{plain}


\bigskip
{\bf \large Contact Information}\\

\begin{tabular}{lcl}
Helmut Prodinger	&& Department of Mathematics, University of Stellenbosch \\ 
warrenham33$@$gmail.com &&7602, Stellenbosch, South Africa
and\\
&&
NITheCS (National Institute for\\&&
Theoretical and Computational Sciences), South Africa.\\
	&& \includegraphics[height=10pt]{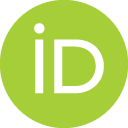}  \url{https://orcid.org/0000-0002-0009-8015}\\
				&& \\
\end{tabular}

\end{document}